\documentclass{amsart}

\usepackage[cp1250]{inputenc}
\usepackage[all]{xy}

\newtheorem{thm}{Theorem}
\newtheorem{prop}[thm]{Proposition}
\newtheorem{crl}[thm]{Corollary}
\newtheorem{lemma}[thm]{Lemma}

\theoremstyle{definition}
\newtheorem*{remark}{Remark}

\DeclareMathOperator{\id}{id}
\DeclareMathOperator{\pd}{d}
\DeclareMathOperator{\Aut}{Aut}
\DeclareMathOperator{\Sl}{sl}
\DeclareMathOperator{\Img}{Im}
\DeclareMathOperator{\hor}{h}
\DeclareMathOperator{\ver}{v}

\newcommand{\sta}{\star}
\newcommand{\stb}{\star_{t}}
\newcommand{\deltb}{\delta_{t}}
\newcommand{\stc}{\star_{a}}
\newcommand{\deltc}{\delta_{a}}

\newcommand{\F}{\mathcal{F}}
\newcommand{\G}{\mathcal{G}}
\newcommand{\R}{\mathbb{R}}

\newcommand{\sld}{\Sl(2)}
\newcommand{\set}[1]{\left\{#1\right\}}
\newcommand{\Lie}[1]{\mathcal{L}_{#1}}
\newcommand{\X}[3]{\mathfrak{X}^{#1}_{#2}\left(#3\right)}
\newcommand{\Om}[3]{\Omega^{#1}_{#2}\left(#3\right)}
\newcommand{\Int}[1]{\iota_{#1}}

\title{A Remark on the Brylinski Conjecture for Orbifolds}
\author{Lukasz Bak}
\address{L. Bak, Institute of Mathematics, Jagiellonian University, ul. Lojasiewicza 6, Cracow, Poland}
\email{lukasz.bak@im.uj.edu.pl}

\author{Andrzej Czarnecki}
\address{A. Czarnecki, Institute of Mathematics, Jagiellonian University, ul. Lojasiewicza 6, Cracow, Poland}
\email{andrzej.czarnecki@im.uj.edu.pl}

\begin{document}

\begin{abstract}
We present reformulation of Mathieu's result on representing cohomology classes of symplectic manifold with symplectically harmonic forms. We apply it to the case of foliated manifolds with transversally symplectic structure and to symplectic orbifolds. We obtain in particular that such representation is always possible for compact K\"{a}hler orbifolds.
\end{abstract}

\maketitle

\section{Introduction}

The goal of this short paper is to prove an orbifold version of the Brylinski conjecture \cite{Bry88}. It concerns the question whether on compact symplectic manifold every cohomology class admits a symplectically harmonic representative. Brylinski proved it in some interesting cases, most notably for compact K\"{a}hler manifolds, but in general the conjecture is not true. The first to present the counterexample was Mathieu (cf. \cite{Mat95}). Moreover, he gave an equivalent condition, in terms of cohomological properties of the symplectic form, for a symplectic manifold (not necessarily compact) to satisfy the conjecture. The same question of the existence of symplectically harmonic representatives can be stated for symplectic orbifolds, and in particular compact K\"{a}hler orbifolds. The conclusion of this paper is the following

\begin{thm} \label{kahler}
Brylinski conjecture holds for compact K\"{a}hler orbifolds.
\end{thm}

The study of homological properties of singular spaces tends to be more complicated than in the nonsingular case. Many properties, which prove very useful in the smooth setting, such as Poincar\'{e} duality or finite dimension of homology groups do not hold in general. These problems have been addressed in different ways. One approach is to refine the homology theory in question to fit this singular spaces. This approach led to definitions of new homologies, like for example Goresky and MacPherson's intersection homology and cohomology (cf. \cite{GorMac80}). For orbifolds, which is of interest to us, another method is applicable.

It is well known (cf. \cite{Mol86}) that the space of leaves of a Riemannian foliation of a compact manifold with compact leaves is an orbifold. In \cite{GirHaeSun83} authors proved the converse, namely that

\begin{thm}
Every orbifold can be realized as the space of leaves of a Riemannian foliation.
\end{thm}

Their proof is written for complex orbifolds and transversally holomorphic foliations but it can be adapted to the real case. This construction, despite its convenience for research of orbifolds, has not been fully exploited. Recently it was used in \cite{WanZaf09} to give a simple proof of Hard Lefschetz Theorem for K\"{a}hler orbifolds, the fact which will be useful for us. Earlier, few others authors took up the similar idea of using foliated Riemannian manifolds and their foliated objects to the study of geometry of the leaf (closure) space (cf. \cite{AleKriLosMic03, JozWol07, PonRadWol09}). The construction allows us to consider ``transverse'' objects on a foliated manifold rather than the corresponding objects on an orbifold. Our next step is to reformulate the problem again, this time using well-known correspondence between ``transverse'' objects on foliated manifolds and ``holonomy invariant'' objects on the corresponding transverse manifolds. We will only sketch this correspondence in a scope necessary for this paper. More general and exhaustive approach can be found for example in \cite{Wol89}. Eventually, we shall see that Theorem \ref{kahler} is a corollary of the following

\begin{thm} \label{main}
Let $M$ be a manifold and $\Gamma$ be a pseudogroup of local diffeomorphisms of $M$. Let $\omega$ be a $\Gamma$-invariant symplectic form on $M$. Then the following conditions are equivalent:
\begin{enumerate}
\item every $\Gamma$-invariant cohomology class has a harmonic representative;
\item for each $k \in \set{0,1,\ldots,n}$ the mapping $L^{k} :\ H^{n-k}_{\Gamma}(M) \to H^{n+k}_{\Gamma}(M)$ is surjective where $L[\xi] = [\omega \wedge \xi]$.
\end{enumerate}
\end{thm}

\section{Invariant forms}

We shall consider a smooth manifold $M$ of dimension $\dim M = m$ together with a pseudogroup of local diffeomorphisms $\Gamma$. From the complex of differential forms $\Om{\ast}{}{M}$ we single out \emph{$\Gamma$-invariant forms}, that is $\xi \in \Om{\ast}{}{M}$ satisfying
\[ \forall \gamma \in \Gamma : \ \xi|_{U} = \gamma^{\ast} \left( \xi|_{\gamma(U)} \right), \]
where $U$ is the domain of $\gamma$. They form a subcomplex (which we shall denote $\Om{\ast}{\Gamma}{M}$) with differential $\pd = \pd|_{\Om{\ast}{\Gamma}{M}}$. Homology of this complex is called \emph{$\Gamma$-invariant cohomology} and denoted $H^{\ast}_{\Gamma}(M)$. In a similar manner we can define \emph{$\Gamma$-invariant multivector fields} as those $X \in \X{\ast}{}{M}$ satisfying
\[ \forall \gamma \in \Gamma : \ X|_{\gamma(U)} = \gamma_{\ast} \left( X|_{U} \right). \]
Again we shall use notation $\X{\ast}{\Gamma}{M}$. For $k < l$ we have the natural nondegenerate pairing $\X{k}{}{M} \times \Om{l}{}{M} \ni (X,\xi) \mapsto \Int{X} \xi \in \Om{l-k}{}{M}$. Direct computations show that this pairing restricts to the pairing $\X{k}{\Gamma}{M} \times \Om{l}{\Gamma}{M} \to \Om{l-k}{\Gamma}{M}$. Using a $\Gamma$-invariant volume form, we obtain an isomorphism

\begin{equation} \label{isom1}
\X{k}{\Gamma}{M} \cong \Om{m-k}{\Gamma}{M}.
\end{equation}

From now on we shall assume that $M$ is of even dimension, $m = 2n$. A $\Gamma$-invariant symplectic form is just a closed, nondegenerate $\omega \in \Om{2}{\Gamma}{M}$. Then the form $\omega^{n}$ is a $\Gamma$-invariant volume form and isomorphism (\ref{isom1}) follows. On the other hand, we have an isomorphism $\X{}{\Gamma}{M} \ni X \mapsto \Int{X} \omega \in \Om{1}{\Gamma}{M}$ which extends to an algebra isomorphism $\X{\ast}{\Gamma}{M} \cong \Om{\ast}{\Gamma}{M}$. Combining these two we obtain an operator
\[ \sta : \  \Om{\ast}{\Gamma}{M} \cong \X{\ast}{\Gamma}{M} \cong \Om{2n-\ast}{\Gamma}{M} \]
associated with the symplectic structure $\omega$.

\begin{remark}
Obviously, operator $\sta$ is an isomorphism. Moreover, it satisfies $\sta^{2} = \id_{\Om{\ast}{\Gamma}{M}}$.
\end{remark}

Now consider a codifferential $\delta$ given, for $\xi \in \Om{k}{\Gamma}{M}$, by
\[ \delta \xi = (-1)^{k} \sta \pd \sta \xi. \]
By the remark above $\delta^{2} = 0$. The central notion in our further studies will be the notion of a \emph{harmonic $\Gamma$-invariant form}: a form $\xi \in \Om{\ast}{\Gamma}{M}$ such that $\pd \xi = 0$ and $\delta \xi = 0$. These forms constitute a subalgebra of $\Om{\ast}{\Gamma}{M}$ which we shall denote by $\Om{\ast}{\Gamma,h}{M}$. The proof of Theorem \ref{main} will follow closely Yan's proof of Mathieu's result (cf. \cite{Yan96}). Let us recall Yan's method. We omit some calculations identical to those in \cite{Yan96}.

We construct a Lie algebra $\sld$-representation on $\Om{\ast}{\Gamma}{M}$. We represent the usual basis of $\sld$ given by
\[ X=\left[ \begin{array}{cc} 0 & 1 \\ 0 & 0 \end{array} \right],\ Y=\left[ \begin{array}{cc} 0 & 0 \\ 1 & 0 \end{array} \right],\ H=\left[ \begin{array}{cc} 1 & 0 \\ 0 & -1 \end{array} \right] \]
as operators on $\Om{\ast}{\Gamma}{M}$ (denoted by the same letters) so that the relations $[X,Y]=H$, $[H,X]=2X$ and $[H,Y]=-2Y$ hold. We set
\begin{align*}
Y \leadsto \Om{k}{\Gamma}{M} \ni \xi & \mapsto \omega \wedge \xi \in \Om{k+2}{\Gamma}{M}, \\
X \leadsto \Om{k}{\Gamma}{M} \ni \xi & \mapsto \sta Y \sta \xi \in \Om{k-2}{\Gamma}{M}, \\
H \leadsto \Om{k}{\Gamma}{M} \ni \xi & \mapsto (n-k) \xi \in \Om{k}{\Gamma}{M}.
\end{align*}

Observe that operator $L$ from the statement of Theorem \ref{main} is induced by $Y$. Eigenvectors of $H$ from the kernel of $X$ are called \emph{primitive}. We say that $\sld$-representation on some vector space $V$ is \emph{of finite $H$-spectrum} iff this space splits into finite sum of eigenspaces of operator $H$. This is obviously the case for our representation. For $\sld$-representations of finite $H$-spectrum the following properties hold (cf. \cite{GriHar78, Yan96})

\begin{prop} \label{props}
Let $V$ be a vector space with $\sld$-representation of finite $H$-spectrum. Let $V_{\lambda}$ be the eigenspace of $H$ of eigenvalue $\lambda$ and $P_{\lambda}$ stand for the set of primitive elements in the eigenspace $V_{\lambda}$. Then
\begin{enumerate}
\item all eigenvalues of $H$ are integers; \item $X|_{V_{k}} : \ V_{k} \to V_{k+2}$, $Y|_{V_{k}} : \ V_{k} \to V_{k-2}$, $k \in \mathbb{Z}$;
\item $X^{k}|V_{-k} : \ V_{-k} \to V_{k}$ and $Y^{k}|V_{k} : \ V_{k} \to V_{-k}$ are isomorphisms for each $k \in \mathbb{N}$;
\item $P_{k} = \set{v \in V_{k} : \ Y^{k+1} v = 0}$;
\item every $V_{k}$ admits a decomposition into the direct sum $V_{k} = P_{k} \oplus \Img Y|_{V_{k+2}}$.
\end{enumerate}
\end{prop}

Applying the above in our case, we obtain following

\begin{crl}
Operator $Y^{k} : \ \Om{n-k}{\Gamma}{M} \to \Om{n+k}{\Gamma}{M}$ is an isomorphism.
\end{crl}

The next step is to prove that the representation on $\Om{\ast}{\Gamma}{M}$ induces an $\sld$-representation on the subspace of $\Gamma$-invariant harmonic forms $\Om{\ast}{\Gamma, h}{M}$, i.e. that $X$, $Y$ and $H$ preserve harmonic forms. This follows from the relations
\[[Y,d] = [X,\delta] = 0,\ [X,d] = - \delta,\ [Y,\delta] = -d,\ [H,d] = [H,\delta] = 0. \]
Therefore we get

\begin{crl}
Operator $Y^{k} : \ \Om{n-k}{\Gamma,h}{M} \to \Om{n+k}{\Gamma,h}{M}$ is an isomorphism.
\end{crl}

\begin{proof} [Proof of Theorem \ref{main}]
We are now ready to prove implication from (1) to (2) of Theorem \ref{main}. The following diagram
\[\xymatrix{
\Om{n-k}{\Gamma,h}{M} \ar[r]_{\cong}^{Y^{k}} \ar[d] & \Om{n+k}{\Gamma,h}{M} \ar[d] \\
H^{n-k}_{\Gamma}(M) \ar[r]^{L^{k}} & H^{n+k}_{\Gamma}(N) }\]
with obvious vertical arrows, is clearly commutative. If we assume surjectivity of vertical arrows, i.e. (1), we obtain surjectivity of lower horizontal arrow. In the proof of the converse, we will use the following lemma, which is a simple corollary of the relation $\delta = -[X,d]$.

\begin{lemma} \label{prim}
A form $\xi \in \Om{\ast}{\Gamma}{M}$ which is closed and primitive must be harmonic.
\end{lemma}

We now assume (2) and proceed with the proof of (1) by induction. Each closed $0$-form is clearly harmonic. It follows from Lemma \ref{prim} that every closed $1$-form is harmonic. Now let $k$ be such that $n - k \geq 2$ and (1) holds for every cohomology class of degree $r < n - k$. Let now $\xi \in \Om{n-k}{\Gamma}{M}$ be a closed form. It suffices to prove that $\xi$ is cohomologous to some harmonic form. By (2) there exists a closed form $\eta \in \Om{n-k-2}{\Gamma}{M}$ and a form $\theta \in \Om{n+k+1}{\Gamma}{M}$ such that $\omega^{k+2} \wedge \eta + \pd \theta = \omega^{k+1} \wedge \xi$. Form $\eta$ is cohomologous to some harmonic form $\bar{\eta} \in \Om{n-k-2}{\Gamma,h}{M}$, so $\eta = \bar{\eta} + \pd \lambda$ for some $\lambda \in \Om{n-k-3}{\Gamma}{M}$. We know that $Y^{k+1} : \ \Om{n-k-1}{\Gamma}{M} \to \Om{n+k+1}{\Gamma}{M}$ is surjective so we can pick $\zeta \in \Om{n-k-1}{\Gamma}{M}$ satysfying $\omega^{k+1} \wedge \zeta = \theta$. Finally, we obtain
\[ \omega^{k+1} \wedge \left[\xi - \pd \left(\zeta + \omega \wedge \lambda\right) - \omega \wedge \bar{\eta}\right] = 0. \]
The form $\bar{\xi} = \xi - \pd \left(\zeta + \omega \wedge \lambda\right)$ is cohomologous to $\xi$. By (4) of Proposition \ref{props} $\bar{\xi} - \omega \wedge \bar{\eta}$ is primitive, therefore by Lemma \ref{prim}, it is harmonic. But $\omega \wedge \bar{\eta}$ is harmonic, and the harmonicity of $\bar{\xi}$ follows.
\end{proof}

\section{Foliations}

Let $M$ be a manifold with a regular foliation $\F$ of dimension $p$ and even codimension $2n$. For a sufficiently small open set $U\subset M$ there is a submersion $p : \ U \to \R^{2n}$ such that foliation $\F$ restricted to $U$ is given by fibers of $p$. Let $\set{U_{i}}$ be an atlas of $M$ consisting of open sets admitting submersions $p_{i} : \ U_{i} \to \R^{2n}$ as above. With this atlas we can associate the transverse manifold $N = \coprod p_{i}(U_{i})$ with holonomy pseudogroup $\Gamma$ generated by the mappings from the associated Haefliger cocycle (cf. \cite{MoeMrc03}).

Consider the basic complex $\Om{\ast}{B}{M,\F}$ consisting of \emph{basic forms} $\xi$ defined by condition $\Int{T}\xi = \Lie{T}\xi = 0$ for each vector field $T$ tangent to the foliation. The cohomology of this complex is also called basic and is denoted by $H^{\ast}_{B}(M,\F)$. Intuitively, basic forms are those which locally are pullbacks of forms from the transverse manifold. This intuition is justified by the following 

\begin{prop} \label{rel}
There is a chain isomorphism $\Xi : \ \Om{\ast}{B}{M,\F}
\stackrel{\cong}{\to} \Om{\ast}{\Gamma}{N}$ satisfying
\[\xi|_{U_{i}} = p_{i}^{\ast}
\left(\Xi(\xi)|_{p_{i}(U_{i})}\right)\] for every $\xi \in
\Om{\ast}{B}{M,\F}$. It is, moreover, an algebra homomorphism.
\end{prop}

Fix a transversely symplectic structure on $M$, that is a basic, closed, nondegenerate $2$-form $\omega$. By Proposition \ref{rel} it induces a $\Gamma$-invariant symplectic structure on $N$, namely $\Xi(\omega)$. Consider the mapping $\X{}{}{M} \ni X \mapsto \Int{X} \omega \in \Om{1}{}{M}$. To make it an isomorphisms we shall restrict it to the \emph{basic vector fields}, defined as follows. Let $\mathcal{D}$ be an arbitrary distribution in $TM$ complementary to $T\F$. The space of basic vector fields $\X{}{B}{M,\F}$ is given by
\[ \X{}{B}{M,\F} = \set{X \in \X{}{}{\mathcal{D}} : \ [X,T] \in \X{}{}{\F}, \forall T \in \X{}{}{\F} } \]
where $\X{}{}{\F}$ stands for vector fields tangent to the foliation and $\X{}{}{\mathcal{D}}$ for those in $\mathcal{D}$. It is now a simple calculation to check that
\[ \X{}{B}{M,\F} \ni X \mapsto \Int{X} \omega \in \Om{1}{B}{M,\F} \]
is a well-defined isomorphism. If we take the subalgebra generated by basic vector fields in $\X{\ast}{}{M}$ and denote it $\X{\ast}{B}{M,\F}$ then this isomorphism extends to an isomorphism $\X{\ast}{B}{M,\F} \cong \Om{\ast}{B}{M,\F}$.

On the other hand
\[ \X{\ast}{B}{M,\F} \ni X \mapsto \Int{X} \omega^{n} \in \Om{2n - \ast}{B}{M,\F} \]
is, again by simple calculations, an isomorphism and we obtain the transversally symplectic star operator
\[ \stb : \ \Om{\ast}{B}{M,\F} \cong \X{\ast}{B}{M,\F} \cong \Om{2n - \ast}{B}{M,\F}. \]
Observe that while both composed isomorphisms depend on the choice of the distribution $\mathcal{D}$, they composition does not. Having the star-operator we proceed as usual, that is we define a codifferential $\deltb \xi = (-1)^{k} \stb \pd \stb$ for $\xi \in \Om{k}{B}{M,\F}$ and consider forms $\xi$ for which $\pd \xi = \deltb \xi = 0$, namely the harmonic forms.

We shall use the chain isomorphism $\Xi$ to apply Theorem \ref{main} in our case. It is possible since the following holds:

\begin{lemma} \label{com}
The following diagram is commutative
\[ \xymatrix{
\Om{\ast}{B}{M,\F} \ar[d]_{\Xi} \ar[r]^{\stb} & \Om{2n - \ast}{B}{M,\F} \ar[d]^{\Xi} \\
\Om{\ast}{\Gamma}{N} \ar[r]^{\sta} & \Om{2n - \ast}{\Gamma}{N}
} \]
\end{lemma}

It assures that $\Xi$ sets a correspondence between basic harmonic forms on $M$ and $\Gamma$-invariant harmonic forms on $N$. Now we can apply Theorem \ref{main} to obtain the following

\begin{thm} \label{main-f}
The following conditions are equivalent:
\begin{enumerate}
\item every basic cohomology class has a harmonic representative;
\item for each $k \in \set{0,1,\ldots,n}$ the mapping $L^{k} :\ H^{n-k}_{B}(M,\F) \to H^{n+k}_{B}(M,\F)$ is surjective.
\end{enumerate}
\end{thm}

\section{Harmonicity with respect to the metric and symplectic structure}

In \cite{Pak08} the author considers the question of existence of harmonic forms when the operator $\sta$ is defined with respect not only to the transversal structure of dimension-one foliation but also some leafwise structure. To avoid confusion we shall denote this new operator by $\stc$. This approach seems proper if we want to define harmonicity in the sense of the whole manifold, not only the transverse structure. We shall now reformulate this problem for the foliations of arbitrary dimension. Recall that $M$ is a manifold with regular foliation $\F$ of dimension $p$ and even codimension $2n$ with transversally symplectic structure $\omega \in \Om{2}{B}{M,\F}$. Assume that leaves of this foliation are orientable or equivalently that $M$ is orientable. Fix a metric $g$ on $M$. Tangent bundle $TM$ splits into direct sum $TM = T\F \oplus T\F^{\bot}$. This splitting induces a bigradation of the space of forms
\[ \Om{\ast}{}{M} = \bigoplus^{p + 2n}_{k = 0} \Om{k}{}{M} = \bigoplus^{p + 2n}_{k = 0} \bigoplus_{r, s \geq 0, r + s = k}\Om{r,s}{}{M,\F} \]
where $\xi \in \Om{r,s}{}{M,\F}$ (is of type $(r,s)$) iff $\Int{X_{1} \wedge \ldots \wedge X_{r+1}} \xi = \Int{Y_{1} \wedge \ldots \wedge Y_{s+1}} \xi = 0$ for $X_{1}, \ldots, X_{r+1} \in \X{}{}{\F}$ and $Y_{1}, \ldots, Y_{s+1} \in \X{\bot}{}{\F}$. Clearly basic $k$-forms are of type $(0,k)$.

Orientability of leaves allows us to choose a leafwise volume form $\chi$ of type $(p,0)$ with respect to $g$. We obtain the volume form $\omega^{n} \wedge \chi$ on $M$, together with an isomorphism
\[ \X{\ast}{}{M} \ni X \mapsto \Int{X} \left(\omega^{n} \wedge \chi\right) \in \Om{p + 2n - \ast}{}{M}. \]
Like before, we define an isomorphism $\X{}{}{M} \to \Om{1}{}{M}$, denoted $\flat$ and given by
\[ \flat (X) = \Int{X^{\ver}} \omega + g(X^{\hor},\cdot) \]
for a decomposition $X = X^{\hor} + X^{\ver} \in \X{}{}{\F} + \X{\bot}{}{\F}$, and then extend it to $\flat : \ \X{\ast}{}{M} \to \Om{\ast}{}{M}.$ The star operator obtained as usual from composition
\begin{equation} \label{st}
\stc : \ \Om{\ast}{}{M} \cong \X{\ast}{}{M} \cong \Om{p + 2n - \ast}{}{M}
\end{equation}
depends on the choice of the metric $g$. We will be interested in the behavior of this operator restricted to basic forms. Consider the bigradation $\X{\ast,\ast}{}{M}$ such that $(r + s)$-vector field $X$ is of type $(r,s)$ iff $\Int{X} \xi = 0$ for each $(r + s)$-form $\xi$ of type other than $(r,s)$. We can now rewrite (\ref{st}) as
\[ \stc : \ \Om{\ast,\ast}{}{M,\F} \cong \X{\ast,\ast}{}{M} \cong \Om{p - \ast, 2n - \ast}{}{M}. \]
For a basic $k$-form $\xi$ we obtain that $\flat^{-1}(\xi)$ is of type $(0,k)$ therefore
\[ \stc \xi = \Int{\flat^{-1}(\xi)} \left( \omega^{n} \wedge \chi \right) = \Int{\flat^{-1}(\xi)} \omega^{n} \wedge \chi = \stb \xi \wedge \chi. \]

If we now consider codifferential defined for $k$-forms by $\deltc \xi = (-1)^{k} \stc \pd \stc \xi$, then for a basic $k$-form we obtain
\begin{align} \label{codiff}
\deltc \xi & = (-1)^{k} \stc \pd \stc \xi = (-1)^{k} \stc \pd \left( \stb \xi \wedge \chi \right) = \stc \left( (-1)^{k} \pd \stb \xi \wedge \chi + \stb \xi \wedge \pd \chi \right) \notag \\
 & = \Int{(-1)^{k} \flat^{-1} (\pd \stb \xi) \wedge \flat^{-1} (\chi)} \left( \omega^{n} \wedge \chi \right) + \Int{\flat^{-1} (\stb \xi) \wedge \flat^{-1} (\pd \chi)} \left( \omega^{n} \wedge \chi \right) \\
  & = \Int{\flat^{-1} (\chi)} \left( \deltb \xi \wedge \chi \right) + \Int{\flat^{-1}(\pd \chi)} \left( \stb^{2} \xi \wedge \chi \right) = (-1)^{p(k-1)} \deltb \xi + \Int{\flat^{-1}(\pd \chi)} \left( \xi \wedge \chi \right) \notag
\end{align}

Observe that with respect to the bigradation $\Om{\ast,\ast}{}{M}$ differential $\pd$ splits into three parts, $\pd = \pd_{1,0} + \pd_{0,1} + \pd_{-1,2}$. In particular $\pd \chi = \pd_{0,1} \chi + \pd_{-1,2} \chi$, for $\pd_{1,0} \chi$ vanishes. The summand $\pd_{-1,2} \chi$ does not vanish in general and therefore $\deltc \xi$ is not basic for a basic form $\chi$. We must therefore redefine the codifferential as $\deltc \xi = (-1)^{k} \stc \pd_{0,1} \stc \xi$. The codifferential is then an adjoint of $\pd_{0,1}$ instead of $\pd$, but on basic forms these two operators coincide. Now if we recall that $\pd_{0,1} \chi = \kappa \wedge \chi$ where $\kappa$ is the mean curvature form for $\F$ and denote $T = \flat^{-1}(\kappa)$ we get
\[ \deltc \xi = (-1)^{p(k-1)} \left(\deltb \xi + \Int{T} \xi \right) . \]
So for basic forms to be preserved by codifferential we must assume that $\kappa$ is basic. Foliations with this property are called \emph{tense}. Recall the operator
\[ Y : \ \Om{\ast}{B}{M,\F} \ni \xi \mapsto \omega \wedge \xi \in \Om{\ast}{B}{M,\F}. \]
From Lemma \ref{com} we get that $[Y, \deltb] = -\pd$. We will compute $[Y, \deltc] \xi$ for a basic $k$-form $\xi$.
\begin{align*}
[Y, \deltc] \xi & = Y \deltc \xi - \deltc Y \xi = (-1)^{p(k-1)} \left( [Y, \deltb] \xi + Y \Int{T} \xi - \Int{T} Y \xi \right) \\
 & = (-1)^{p(k-1)+1} \left( \Int{T} \omega \wedge \xi + \pd \xi \right) .
\end{align*}
This means that to obtain the invariance of basic forms by $Y$, we need to assume that $\Int{T} \omega = 0$. But then nondegeneracy of $\omega$ implies vanishing of the mean curvature. In this case basic forms harmonic with respect to $\stc$ and basic forms harmonic with respect to $\stb$ are the same, because $\stc = \stb$. So we can conclude these considerations with

\begin{thm} \label{main-ft}
For the foliation $\F$ with vanishing mean curvature form $\kappa$ the following conditions are equivalent:
\begin{enumerate}
\item every basic cohomology class has a harmonic representative with respect to the operator $\stc$;
\item for each $k \in \set{0,1,\ldots,n}$ the mapping $L^{k} :\ H^{n-k}_{B}(M,\F) \to H^{n+k}_{B}(M,\F)$ is surjective.
\end{enumerate}
\end{thm}

By a well-known technique, we can deform the metric along the leaves to change the mean curvature form inside its cohomology class. In particular, if we can find a metric on the manifold $M$ for which the mean curvature form is exact, then we can find another metric with vanishing mean curvature. Especially, for a tense, transversally oriented Riemannian foliation on a closed oriented manifold, exactness of mean curvature form is equivalent to the minimalizability of leaves (cf. \cite{Ton97}). Due to \cite{Dom98} we can omit the tenseness condition. Recall that we consider foliations with orientable leaves. We can state the following

\begin{crl}
For a transversally symplectic, minimalizable Riemannian foliation on a closed manifold there exist a metric, for which conditions (1) and (2) of Theorem \ref{main-ft} are equivalent.
\end{crl}

In \cite{Mas92} Masa proves that for a transversally orientable Riemannian foliation on closed, orientable manifold minimalizability of leaves is equivalent to nontriviality of highest rank basic cohomology. Applying it to our case, we obtain

\begin{crl}
For a transversally symplectic Riemannian foliation of codimension $2n$, on a closed manifold, such that $H^{2n}_{B} (M, \F) \neq 0$ there exist a metric, for which conditions (1) and (2) of Theorem \ref{main-ft} are equivalent.
\end{crl}

\section{Orbifolds}

The notion of the orbifold is a generalization of the notion of manifold. While a manifold locally looks like an Euclidean space, an orbifold locally looks like a quotient space of an Euclidean space under an action of a finite group. This notion includes manifolds, but also manifolds with boundary and even more, like manifolds with corners. The first definition of an orbifold was given by Satake in \cite{Sat56} under the name of V-manifold. Numerous other definitions were given later, among them was the one given by Thurston in \cite{Thu02} or by Moerdijk and Mrcun in \cite{MoeMrc03}. In this note we follow the definitions of Chen and Ruan (\cite{CheRua02}) using local uniformizing systems. We recall them in brief.

An \emph{$n$-dimensional local uniformizing system} over some open set $U \subset M$ is a triple $(V, \G, \pi)$ consisting of an open set $V \subset \R^{n}$, a finite group $\G$ acting smoothly and effectively on $V$ and a continous map $\pi : \ V \to U$, invariant by the action of $\G$, inducing a homeomorphism between $V/G$ and $U$. An open covering $\set{U_{i}}$ of a Hausdorff, second-countable topological space together with a family of local uniformizing systems $\set{(V_{i}, \G_{i}, \pi_{i})}$ satisfying some compatibility conditions is called an orbifold structure on $X$.

A surjective, continous map $p : \  E \to U$ is a local trivialization of an orbifold bundle of rank $k$ if there are local uniformizing systems $(V, \G, \pi)$ for $U$ and $(V \times \R^{k}, \G, \tilde{\pi})$ for $E$ and a smooth map $\rho : \ V \times \G \to \Aut(\R^{k})$ that satisfies the following group condition:
\[ \forall x \in V, \forall g,h \in \G : \ \rho(gx, h) \rho(x, g) = \rho(x, hg), \]
The action of $\G$ on $V \times \R^{k}$ is given by $g(x,v) = (gx, \rho(x,g) v)$ and the following diagram is commutative:
\[ \xymatrix{
V \times \R^{k} \ar[r] \ar[d]_{\tilde{\pi}} & V \ar[d]^{\pi} \\
E \ar[r]^{p} & U.
} \]
Again, the family of local trivialisations over some open covering of $X$ satisfying some compatibility conditions form an uniform bundle.

We can consider any orbifold $X$ with an orbifold structure $\set{(V_{i}, \G_{i}, \pi_{i})}$ over some open covering $\mathcal{U} = \set{U_{i}}$. By passing to a refined covering we can assume that each $V_{i}$ is an open subset of $\R^{n}$. Hence we obtain an uniformizing bundle system of rank $n$ over $U_{i}$ by taking $(TV_{i}, \G_{i}, \tilde{\pi}_{i})$ and $\rho_{i}(v,g) = \pd_{v}g$. Any two such systems are equivalent, and so they form an orbifold bundle $TX$, which plays role of a tangent bundle. In a similar way we can obtain the cotangent bundle $T^{\ast}X$ and in general any tensor bundle. We can consider sections of the bundle $E$, namely the lifts of the mappings $s : \ X \to E$ which in local uniformizing bundle systems are right inverse of the natural projection $V \times \R^{k} \to V$. In case of tensor bundles those will be simply tensor fields on the orbifold $X$, in particular sections of bundle $\bigwedge^{k}T^{\ast}X$ are just $k$-forms, and the space of such sections we will denote as $\Om{k}{}{X}$. Now, the operation $\pd : \ \Om{k}{}{X} \to \Om{k+1}{}{X}$ which in local uniformizing system $(V_{i}, \G_{i}, \pi_{i})$ is given by $\pd : \ \Om{k}{}{V_{i}} \to \Om{k+1}{}{V_{i}}$ is well-defined and satisfies $\pd^{2} = 0$. This is how cohomology groups of the orbifold $X$ are constructed. Observe that the wedge product, which is pointwise, is well-defined.

We can also consider a section of the bundle of positive, symmetric tensors of type $(0,2)$ on the orbifold $X$. Any such section plays a role of Riemannian metric on $X$, and any orbifold $X$ admits such metrics. If we choose a Riemannian metric on $X$, we can consider the bundle of orthonormal frames over $X$. For any uniformizing system $(V,\G,\pi)$ action of $\G$ lifts in natural way to the free action on $LV$, the bundle of orthonormal frames over $V$. The induced action of $O(n)$ on the manifold $LV/\G$ gives rise to a foliation such that corresponding space of leaves has a natural orbifold structure, isomorphic to $X$. For details of this construction see \cite{MoeMrc03}. Basic forms of this foliation are exactly the forms on the orbifold $X$. In particular if we fix a symplectic structure on $X$, that is a nondegenerate, closed $2$-form on $X$, we obtain a basic symplectic form on the associated foliated manifold. Other constructions from Section 2 can be carried out on $X$ as well, and as a result of Theorem \ref{main-f} we obtain

\begin{thm} \label{orbmain}
For a symplectic orbifold $X$ of dimension $n$ the following conditions are equivalent:
\begin{enumerate}
\item every cohomology class has a harmonic representative;
\item for each $k \in \set{0,1,\ldots,n}$ the mapping $L^{k} :\ H^{n-k}(X) \to H^{n+k}(X)$ is surjective.
\end{enumerate}
\end{thm}

This result may be obtained in a more direct way, by finding some correspondence of the form considered in Proposition \ref{rel}. For a fixed orbifold structure $\set{(V_{i}, \G_{i}, \pi_{i})}$ on $X$, such that the associated covering $\mathcal{U}$ is closed under intersections, consider a manifold $N = \coprod_{i} V_{i}$. Let $\tilde{s}$ be a lift of $s : \ X \to \bigwedge^{k}T^{\ast}X$, that is a $k$-form. It corresponds to the form
\[ \sum_{i} \left. \tilde{s} \right|_{V_{i}} \in \bigoplus_{i} \Om{k}{}{V_{i}} = \Om{k}{}{N}. \]
Now from the definitions above we obtain

\begin{prop} \label{orbisom}
The mapping \[\Om{\ast}{}{X} \ni s \to \sum_{i} \left. \tilde{s} \right|_{V_{i}} \in \Om{k}{}{N} \] is a chain isomorphism onto a subcomplex $\Om{k}{\Gamma}{N}$ consisting of forms invariant by pseudogroup $\Gamma$ generated by injections $(V_{i}, \G_{i}, \pi_{i}) \to (V_{j}, \G_{j}, \pi_{j})$ and the action of each $\G_{i}$ on $V_{i}$.
\end{prop}

Proposition \ref{orbisom} together with Theorem \ref{main} now imply Theorem \ref{orbmain}.

\begin{remark}
For compact K\"{a}hler orbifolds the Hard Lefschetz Theorem, which implies (2) in the Theorem \ref{orbmain} was established in \cite{GirHaeSun83}. Hence Theorem \ref{kahler} follows.
\end{remark}

\bibliography{Artykul}
\bibliographystyle{plain}

\end{document}